\newif\ifJ

\documentclass[10pt]{article}
\usepackage[a4paper]{geometry}

\ifJ
\usepackage{amsmath,amsfonts,amssymb}
\else
\usepackage{amsmath,amsfonts,amsthm,amssymb}
\fi
\usepackage{mathrsfs}
\usepackage[breaklinks,bookmarks=false]{hyperref}
\hypersetup{colorlinks, linkcolor=blue, citecolor=blue,
urlcolor=blue, plainpages=false, pdfwindowui=false,
pdfstartview={FitH}}
\usepackage{xcolor}
\usepackage{graphicx} 
\usepackage{mathtools}

\ifJ
\else
\usepackage{caption}
\fi

\usepackage{subfig}
\usepackage{tikz}

\usepackage{pgfplots}
\usepackage{pgfplotstable}
\usepgfplotslibrary{groupplots}
\pgfplotsset{compat=1.18}
\pgfplotsset{compat=newest}
\usepackage{bigints}

\pgfplotscreateplotcyclelist{my markers}{
    mark=+\\
    mark=x\\
    mark=|\\
    mark=asterisk\\
    mark=diamond\\
    mark=triangle\\
}

\pgfplotsset{every axis legend/.append style=
{
font=\small,%
xticklabel style={font=\small},%
xlabel style={font=\small},%
yticklabel style={font=\small}
},%
ylabel style={font=\small},%
legend style={font=\small},%
title style={font=\small},
scale only axis,
grid=both,
grid style={line width=.2pt, draw=gray!10},
major grid style={line width=.4pt, draw=gray!50},
}

\pgfplotsset{every axis plot/.append style=
{
cycle list name=my markers,
line width=0.8pt,
mark size=3pt,
}
}

\ifJ
\else
\usepackage{cleveref}
\fi
\usepackage{float}

\ifJ
\usepackage[disable]{todonotes}
\else
\usepackage[textsize=scriptsize]{todonotes}
\fi

\ifJ
\else

\fi

\newcommand{\dx}{\mathrm{d}x}
\newcommand{\R}{\mathbb{R}}
\newcommand{\N}{\mathbb{N}}

\newcommand{\norm}[1]{\left\lVert#1\right\rVert}

\definecolor{softblue}{RGB}{76,114,176}
\definecolor{softorange}{RGB}{221,132,82}
\definecolor{softgreen}{RGB}{85,168,104}
\definecolor{softred}{RGB}{196,78,82}
\definecolor{softpurple}{RGB}{129,114,179}
\definecolor{softbrown}{RGB}{147,120,96}
\definecolor{softpink}{RGB}{218,139,195}
\definecolor{softgray}{RGB}{140,140,140}

\ifJ

\title{The Effect of Quadrature on the Convergence of Policy Iteration for Hamilton--Jacobi--Bellman Equations\thanks{I.~Smears and H.~Wells were supported by the Engineering and Physical Sciences Research Council [grant number EP/Y008758/1]}}
\titlerunning{Effect of quadrature on policy iteration}
\author{T.~Hall \and I.~Smears \and E.~S\"uli \and H.~Wells}

\institute{T.~Hall, I.~Smears \& H.~Wells \at
              Department of Mathematics, University College London, London, WC1E 6BT.
              \email{(tom.hall.21,i.smears,h.wells)@ucl.ac.uk; }
\and E.~S\"uli \at
    Mathematical Institute, University of Oxford, Oxford, OX2 6GG.
    \email{suli@maths.ox.ac.uk}
}

\date{Received: date / Accepted: date}

\else
\title{The Effect of Quadrature on the Convergence of Policy Iteration for Hamilton--Jacobi--Bellman Equations}

\author{Thomas Hall, Iain Smears, Endre S\"uli, Harry Wells}
\date{June 23, 2026}
\fi

\begin{document}

\maketitle

\begin{abstract}
    Modern finite element libraries allow users to express partial differential equations directly in variational form, with the added convenience of automatic quadrature selection.
    In the context of Hamilton--Jacobi--Bellman (HJB) equations, automatic quadrature selection can result in nonmatching quadratures between different terms that may lead to loss of convergence of the policy iteration, which is otherwise expected from theory to converge superlinearly.
    The simple remedy of enforcing matching quadrature recovers the expected superlinear convergence. 
\end{abstract}
\section{Introduction}
Ever since the seminal works of Bellman~\cite{Bellman1957} and Howard~\cite{Howard1960}, policy iteration remains a mainstay algorithm for solving Hamilton--Jacobi--Bellman (HJB) equations of optimal control.
It has long been known that it can be viewed as a Newton-type method in function spaces~\cite{PutermanBrumelle1979}.
In many cases the nonlinearity of the problem is nonsmooth, in which case policy iteration can be analyzed as a semismooth Newton method, and the theory shows its asymptotic superlinear convergence rate for a wide range of problem settings and discretization methods~\cite{BokanowskiMarosoZidani2009,SmearsSuli2014}.
Based on ample computational evidence over many decades, the algorithm is generally considered to be robust with respect to various sources of error, such as inexact linear algebra and quadrature.

Recently, we have been asked by some colleagues to explain why their implementations of policy iteration in modern finite element software libraries fail to converge.
Modern finite element software libraries, such as Firedrake~\cite{FiredrakeUserManual} and NGSolve~\cite{NGSolve}, allow the user to express the equations to be solved in variational form, with the library automatically selecting quadratures for the terms of the equation.
As we demonstrate below, loss of convergence can occur if the library selects nonmatching quadrature schemes for different terms of the equation.
This issue appears not to have been observed in older software implementations, where the default choice was to use a single quadrature scheme for all terms in the problem.
Once aware of the issue, the remedy is simply to enforce matching quadrature for the various terms in the equation.
\section{Effect of quadrature on policy iteration}
\label{sec:quadrature}

To present the issue in a simple context, consider the following second-order semilinear HJB equation
\begin{equation}\label{eq:model}
\begin{aligned}
     F[u]\coloneqq - \Delta u + H(x, \nabla u) &= 0 && \text{in } \Omega, \\
     u &= 0 && \text{on }\partial \Omega,
 \end{aligned}
\end{equation}
where $\Omega \subset \R^d$ is a bounded Lipschitz domain, and the Hamiltonian $H: \Omega \times \R^d \to \R$ is a real-valued function defined by
\begin{equation}
\label{eq:hamiltonian}
    H(x, p) \coloneqq \max_{\alpha \in \mathcal{A}}\{b(x, \alpha) \cdot p - f(x, \alpha) \} \qquad \forall (x, p) \in \Omega \times \R^d,
\end{equation}
where $\mathcal{A}$ denotes the set of controls of the underlying control problem, and $b$ and $f$ denote respectively the control-dependent drift and running costs.
In the following we assume only that~$\mathcal{A}$ is a compact metric space, and that $(b,f) \in C( \overline{\Omega} \times \mathcal{A}; \R^d\times \R)$.
In particular, under the above assumptions, the Hamiltonian~$H$ is Lipschitz continuous w.r.t.\ the gradient variable.

\subsection{Policy iteration}
Given an initial guess~$u_0$, the policy iteration algorithm is to compute a sequence of approximations~$u_k$ of~$u$, where, for each $k\in \N$, $k \geq 1$, the approximation~$u_k$ solves the linearized problems 
\begin{equation}\label{eq:policy_iteration}
    \begin{aligned}
-\Delta u_k + b_k\cdot \nabla u_k&=f_k && \text{in }\Omega,
\end{aligned}
\end{equation}
along with $u_k=0$ on $\partial \Omega$, where the coefficients $b_k\coloneqq b(\cdot,\alpha_k(\cdot))$ and $f_k=f(\cdot,\alpha_k(\cdot))$ for some (Lebesgue measurable) policy $\alpha_k\colon\Omega\to \mathcal{A}$ that is chosen such that $\alpha_k(x)$ is a maximizer of $b(x,\alpha)\cdot \nabla u_{k-1}(x) - f(x,\alpha)$ over all~$\alpha \in \mathcal{A}$, for a.e.\ $x\in \Omega$. 

Under the above hypotheses, it is well-known that the algorithm produces a monotonically globally convergent sequence in~$H^1_0(\Omega)$, and moreover the convergence is asymptotically superlinear, i.e.\ the ratio of successive errors $\norm{u-u_{k}}_{H^1(\Omega)} / \norm{u-u_{k-1}}_{H^1(\Omega)}$ converges to~$0$ as $k\to\infty$, unless the exact solution is attained in finitely many steps.
The proof of these facts is easily obtained by minor adaptation of the analysis for second-order fully nonlinear HJB operators from~\cite[Theorem~13]{SmearsSuli2014} to the easier setting of semilinear problems.

\subsection{Discretization with exact integration}

In practical computations, the HJB equation~\eqref{eq:model} and corresponding linearized problems~\eqref{eq:policy_iteration} are solved in discretized form, for instance by a finite element method.
Typically, the coefficients~$(b_k,f_k)$ are then only evaluated at the quadrature nodes that are used to approximate the integral terms of the corresponding finite element scheme.
To illustrate the central issue in the simplest possible setting, we consider below a standard FEM discretization of~\eqref{eq:model}.

Let~$\mathcal{T}_h$ be a conforming simplicial mesh on~$\Omega$, which we assume to be a Lipschitz polytope, and take~$V_h$ to be the $H^1_0(\Omega)$-conforming finite element space consisting of piecewise polynomials of degree at most $p \geq 1$ on $\mathcal{T}_h$. 
The numerical scheme is to find $u_h\in V_h$ that solves 
\begin{equation}
\label{eq:weak}
    \big\langle F[u_h], v_h \big\rangle \coloneqq \int_{\Omega} \left( \nabla u_h {\cdot} \nabla v_h  +  H(x, \nabla u_h) v_h \right)\dx =0 \qquad \forall v_h \in V_h,
\end{equation}
where $F\colon H^1_0(\Omega)\to H^{-1}(\Omega)$ is defined in~\eqref{eq:model} above and $\langle \cdot , \cdot \rangle$ denotes the pairing between~$H^{-1}(\Omega)$ and~$H^1_0(\Omega)$. 
The policy iteration assuming exact integration is defined as follows. 
For $v_h \in V_h$, we denote by $\Lambda[v_h]$ the set of all Lebesgue measurable optimal feedback controls w.r.t. $v_h$, that is, all Lebesgue measurable functions $\alpha: \Omega \to \mathcal{A}$ such that $\alpha(x) \in \mathrm{arg}\max_{\alpha \in \mathcal{A}}\big\{b(x, \alpha) \cdot \nabla v_h(x) - f(x, \alpha) \big\}$ for a.e. $x \in \Omega$. 
Suppose that $u_{h, k-1} \in V_h$, $k \geq 1$, is a current iterate, and define the bilinear form
\begin{equation}
\label{eq:bilinear}
     B_k(w_h, v_h) \coloneqq \int_{\Omega} \left(\nabla w_h \cdot \nabla v_h  +  b(x, \alpha_k(x)) {\cdot} \nabla w_h \,v_h\right) \dx \quad \forall w_h,\, v_h\in V_h,
\end{equation}
where $\alpha_k$ is chosen from $\Lambda[u_{h,k-1}]$. 
The \textbf{nonresidual form} of the algorithm is to find $u_{h,k}\in V_h$ that solves
\begin{equation}
\label{eq:policynonres}
    B_k(u_{h, k}, v_h) = \int_{\Omega} f(x, \alpha_k(x)) v_h \dx \qquad \forall v_h \in V_h.
\end{equation}
%
The \textbf{residual form} of the algorithm defines $u_{h, k} \coloneqq u_{h, k-1} + \delta_{h, k}$, where $\delta_{h, k}\in V_h$ solves
\begin{equation}
\label{eq:policyres}
    B_k(\delta_{h, k}, v_h) = - \langle F(u_{h,k-1}), v_h \rangle \qquad \forall v_h \in V_h.
\end{equation}
Assuming exact integration and exact arithmetic, the forms~\eqref{eq:policynonres} and \eqref{eq:policyres} are equivalent and define the same sequence of iterates.
Indeed, the equivalence of~\eqref{eq:policynonres} and~\eqref{eq:policyres} is a simple consequence of $H(x, \nabla u_{h,k-1}) = b(x, \alpha_k(x)) \cdot \nabla u_{h,k-1} - f(x, \alpha_k(x))$, which follows from the choice of $\alpha_k \in \Lambda[u_{h,k-1}]$.
Without attempting to go into detail, the local superlinear convergence of the iterates defined by~\eqref{eq:policynonres} and, equivalently, by~\eqref{eq:policyres} above, can be shown under the assumptions already stated above and the additional assumption that the discrete linearized problems are uniformly stable.

\subsection{Quadrature}
In practice, the integrals in \eqref{eq:policynonres} and \eqref{eq:policyres} must be approximated by numerical quadrature.
If the nonresidual form~\eqref{eq:policynonres} is used as the starting point for an implementation, then modern software libraries with automatic quadrature selection typically approximate~\eqref{eq:policynonres} with an update equation of the form
\begin{multline}
\label{eq:quadnonres}
    \sum_{K \in \mathcal{T}_h} \sum_{i=1}^{N_K} \omega_{i,K} \big[ \nabla u_{h, k}(x_{i, K}) \cdot \nabla v_h(x_{i, K}) \\
    + b(x_{i, K}, \alpha_k(x_{i, K})) \cdot \nabla u_{h, k}(x_{i, K}) v_h(x_{i, K}) \big] \\
    = \sum_{K \in \mathcal{T}_h} \sum_{i=1}^{\widetilde{N}_K} \widetilde{\omega}_{i, K} f(\widetilde{x}_{i, K}, \alpha_k(\widetilde{x}_{i, K})) v_h(\widetilde{x}_{i, K})
    \quad \forall v_h\in V_h,
\end{multline}
where, for each $K \in \mathcal{T}_h$, the sets of pairs $\{(x_{i,K},\omega_{i,K})\}_{i=1}^{N_K}$ and $\{(\widetilde{x}_{i,K},\widetilde{\omega}_{i,K})\}_{i=1}^{\widetilde{N}_K}$ are quadrature nodes and weights used to approximate the integrals over each mesh-element in the left-hand side; respectively the right-hand side, of~\eqref{eq:policynonres}.
In practice, the possible difference in choices of nodes and weights between LHS and RHS may arise from the separate compilation and assembly of the stiffness matrix and right-hand side vector.
In the following, we will call the choice of quadrature \textbf{matching} if the $\{(x_{i,K},\omega_{i,K})\}_{i=1}^{N_K}$ and $\{(\widetilde{x}_{i,K},\widetilde{\omega}_{i,K})\}_{i=1}^{\widetilde{N}_K}$ coincide, otherwise we will call it \textbf{nonmatching}.
If the implementation derives from the residual form~\eqref{eq:policyres}, then quadrature approximation leads to an equation for the update $\delta_{h,k}=u_{h,k}-u_{h,k-1}$ of the form
\begin{multline}
\label{eq:quadres}
    \sum_{K \in \mathcal{T}_h} \sum_{i=1}^{N_K} \omega_{i,K} \big[ \nabla \delta_{h, k}(x_{i, K}) \cdot \nabla v_h(x_{i, K}) \\
    + b(x_{i, K}, \alpha_k(x_{i, K})) \cdot \nabla \delta_{h, k}(x_{i, K}) v_h(x_{i, K}) \big] \\
    = - \sum_{K \in \mathcal{T}_h} \sum_{i=1}^{\widetilde{N}_K} \widetilde{\omega}_{i,K} \big[ \nabla u_{h, k-1}(\widetilde{x}_{i, K}) \cdot \nabla v_h(\widetilde{x}_{i, K}) \\
    + H(\widetilde{x}_{i, K}, \nabla u_{h, k-1}(\widetilde{x}_{i, K})) v_h(\widetilde{x}_{i, K}) \big],
\end{multline}
for all $v_h \in V_h$.
Observe that if the quadrature is matching, then \eqref{eq:quadnonres} and~\eqref{eq:quadres} are equivalent, but for nonmatching quadratures the equivalence is no longer generally valid and the two forms appear to define different sequences of iterates; this will be confirmed in the experiments below.
\section{Numerical experiment}
\label{sec:numerics}
We present below the results for a test problem implemented in NGSolve~\cite{NGSolve}.
Although we do not report here the details, we also performed analogous experiments in Firedrake~\cite{FiredrakeUserManual} where we observed the same behaviour. The code and data used to generate the numerical results presented here are available at~\cite{HJBQuadratureCode}.
\begin{figure}[!b]
\centering
\begin{tikzpicture}
\begin{groupplot}[
    group style={group size=2 by 1, horizontal sep=0.2cm, y descriptions at=edge left},
    width=0.43\textwidth,
    height=0.34\textwidth,
    scale only axis,
    xlabel={Iteration},
    ylabel={Discrete residual norm},
    legend style={font=\tiny},
    ymode=log,
    xmin=-0.2, xmax=8.2,
    ymin=1e-16, ymax=1e1,
]

\nextgroupplot[title={Nonresidual formulation}, legend pos=south west,cycle list name=my markers]
\addplot+[color = softblue, ] table [x=iter, y=dof_100, col sep=comma] {newton_nr_automatic_2panel_p1.csv};
\addplot+[color = softorange,] table [x=iter, y=dof_324, col sep=comma] {newton_nr_automatic_2panel_p1.csv};
\addplot+[color = softgreen, ] table [x=iter, y=dof_1156, col sep=comma] {newton_nr_automatic_2panel_p1.csv};
\addplot+[color = softred, ] table [x=iter, y=dof_4356, col sep=comma] {newton_nr_automatic_2panel_p1.csv};
\addplot+[color = softpurple,] table [x=iter, y=dof_16900, col sep=comma] {newton_nr_automatic_2panel_p1.csv};
\addplot+[color = softbrown, ] table [x=iter, y=dof_66564, col sep=comma] {newton_nr_automatic_2panel_p1.csv};
\legend{100 dofs, 324 dofs, 1156 dofs, 4356 dofs, 16900 dofs, 66564 dofs}

\nextgroupplot[title={Residual formulation}, legend pos=south west,cycle list name=my markers]
\addplot+[color = softblue, ] table [x=iter, y=dof_100, col sep=comma] {newton_r_automatic_2panel_p1.csv};
\addplot+[color = softorange, ] table [x=iter, y=dof_324, col sep=comma] {newton_r_automatic_2panel_p1.csv};
\addplot+[color = softgreen, ] table [x=iter, y=dof_1156, col sep=comma] {newton_r_automatic_2panel_p1.csv};
\addplot+[color = softred, ] table [x=iter, y=dof_4356, col sep=comma] {newton_r_automatic_2panel_p1.csv};
\addplot+[color = softpurple, ] table [x=iter, y=dof_16900, col sep=comma] {newton_r_automatic_2panel_p1.csv};
\addplot+[color = softbrown,] table [x=iter, y=dof_66564, col sep=comma] {newton_r_automatic_2panel_p1.csv};
\legend{100 dofs, 324 dofs, 1156 dofs, 4356 dofs, 16900 dofs, 66564 dofs}

\end{groupplot}
\end{tikzpicture}

\begin{tikzpicture}
\begin{groupplot}[
    group style={group size=2 by 1, horizontal sep=0.2cm, y descriptions at=edge left},
    width=0.43\textwidth,
    height=0.34\textwidth,
    scale only axis,
    xlabel={Iteration},
    ylabel={Discrete residual norm},
    title style={font=\footnotesize, yshift = -2mm},
    legend style={font=\tiny},
    ymode=log,
    xmin=-0.2, xmax=8.2,
    ymin=1e-16, ymax=1e1,
]

\nextgroupplot[title={Nonresidual formulation}, legend pos=north east,cycle list name=my markers]
\addplot+[color = softblue ] table [x=iter, y=dof_100, col sep=comma] {newton_nr_prescribed_2panel_p1_quad_m4_rhs4.csv};
\addplot+[color = softorange] table [x=iter, y=dof_324, col sep=comma] {newton_nr_prescribed_2panel_p1_quad_m4_rhs4.csv};
\addplot+[color = softgreen ] table [x=iter, y=dof_1156, col sep=comma] {newton_nr_prescribed_2panel_p1_quad_m4_rhs4.csv};
\addplot+[color = softred] table [x=iter, y=dof_4356, col sep=comma] {newton_nr_prescribed_2panel_p1_quad_m4_rhs4.csv};
\addplot+[color = softpurple] table [x=iter, y=dof_16900, col sep=comma] {newton_nr_prescribed_2panel_p1_quad_m4_rhs4.csv};
\addplot+[color = softbrown] table [x=iter, y=dof_66564, col sep=comma] {newton_nr_prescribed_2panel_p1_quad_m4_rhs4.csv};
\legend{100 dofs, 324 dofs, 1156 dofs, 4356 dofs, 16900 dofs, 66564 dofs}

\nextgroupplot[title={Residual formulation}, legend pos=north east,cycle list name=my markers]
\addplot+[color = softblue, ] table [x=iter, y=dof_100, col sep=comma] {newton_r_prescribed_2panel_p1_quad_m4_rhs4.csv};
\addplot+[color = softorange, ] table [x=iter, y=dof_324, col sep=comma] {newton_r_prescribed_2panel_p1_quad_m4_rhs4.csv};
\addplot+[color = softgreen, ] table [x=iter, y=dof_1156, col sep=comma] {newton_r_prescribed_2panel_p1_quad_m4_rhs4.csv};
\addplot+[color = softred,] table [x=iter, y=dof_4356, col sep=comma] {newton_r_prescribed_2panel_p1_quad_m4_rhs4.csv};
\addplot+[color = softpurple,] table [x=iter, y=dof_16900, col sep=comma] {newton_r_prescribed_2panel_p1_quad_m4_rhs4.csv};
\addplot+[color = softbrown, ] table [x=iter, y=dof_66564, col sep=comma] {newton_r_prescribed_2panel_p1_quad_m4_rhs4.csv};
\legend{100 dofs, 324 dofs, 1156 dofs, 4356 dofs, 16900 dofs, 66564 dofs}

\end{groupplot}
\end{tikzpicture}
\caption{Convergence histories of the discrete $H^{-1}$-norms of residuals under default automatic quadrature (top) and prescribed matching quadrature (bottom) for the nonresidual form (left) and the residual form (right).}
\label{fig:newton}
\end{figure}
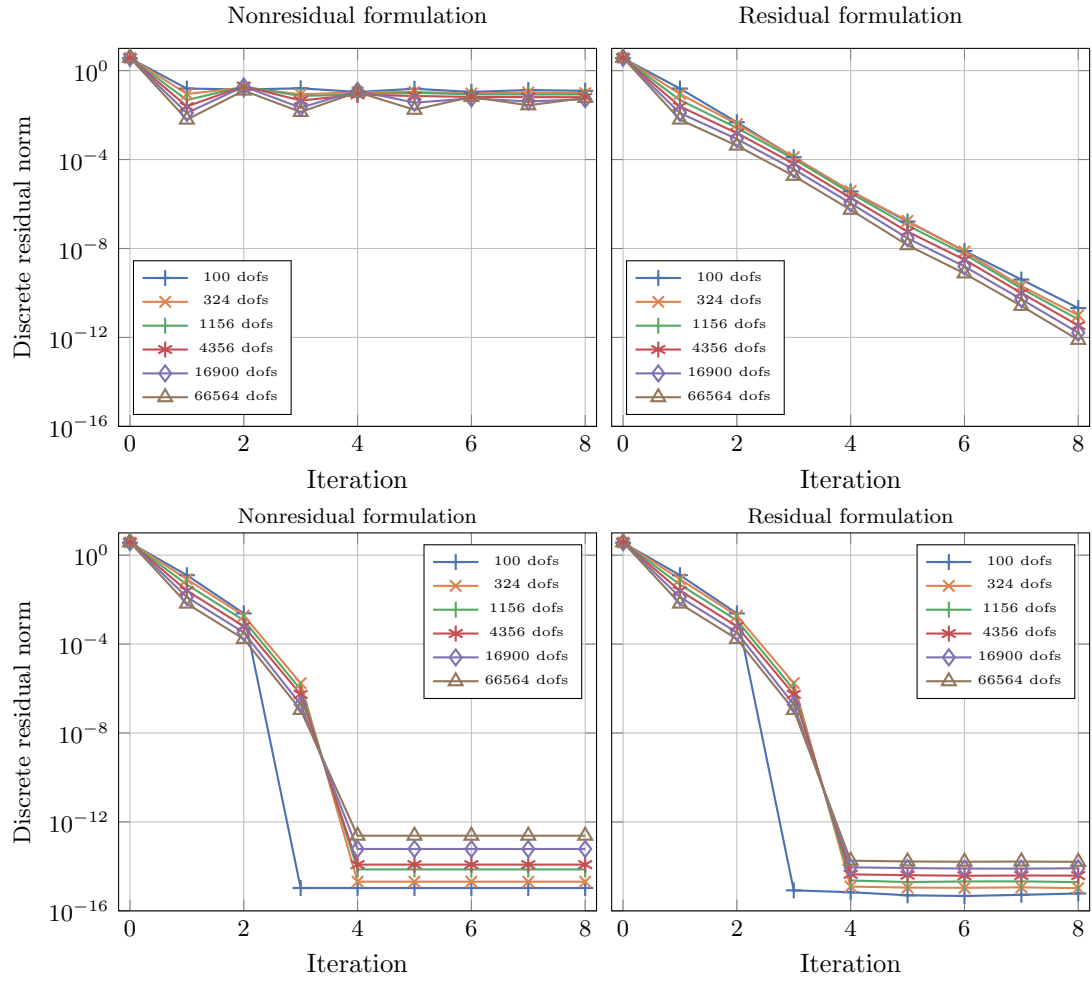
The test problem setup is as follows. Let $\Omega = (-\pi, \pi)^2$, and consider the HJB equation~\eqref{eq:model} where $\mathcal{A} = \{0, 1\}$, and $b\big(x, \alpha\big) = \alpha(\cos x_1 , 0)$ for $x=(x_1,x_2)\in \Omega$ and $\alpha\in \mathcal{A}$.
Let $f(x,\alpha)\coloneqq \alpha \left(\cos^2 x_1\sin x_2-g(x)  \right) + 2\sin x_1 \sin x_2 + \max(g(x),0)$ for all $x=(x_1,x_2)\in \Omega$ and $\alpha\in \mathcal{A}$, where the function $g(x)=10^{-4}\sin(4x_1)\sin(4x_2)$; with this choice of data the exact solution is $u(x)=\sin x_1 \sin x_2$ and the optimal controls for the exact solution jump between the two values in~$\mathcal{A}$ in multiple subregions of~$\Omega$.
This is an example of a nondifferentiable Hamiltonian, where the optimal control that achieves the maximum in~\eqref{eq:hamiltonian} is in general discontinuous with respect to the gradient variable.
We use the FEM defined in~\eqref{eq:weak} with piecewise affine elements ($p=1$) on structured triangular meshes that are obtained from a uniform Cartesian grid by subdividing each square along the bottom-left to top-right diagonal. 

We compare the results obtained with the default automatic quadrature selection and those with a manually prescribed quadrature selection. 
The implementations for default and prescribed quadratures differ only by a single line of code.
The prescribed quadrature for the results reported here was of order 4, although similar results are found for all orders from 2 to 12 that we tested.
Figure~\ref{fig:newton} shows that, for the default automatic quadrature selection, there is no convergence at all for the nonresidual form, whereas the residual formulation converges albeit only at a linear rate.
The observed loss of equivalence between residual and nonresidual forms in the results demonstrates that there is a mismatch in the quadrature choices between the terms in the bilinear forms and right-hand sides.
The linear convergence of the residual form~\eqref{eq:quadres} resembles that of inexact Newton methods, where it is well-known that one can allow some error in the differential of the nonlinear operator whilst still retaining linear convergence.
Figure~\ref{fig:newton} also shows that, for matching quadrature, both residual and nonresidual forms of policy iteration converge with superlinear rate, which corresponds to the rate expected from theory.

\section{Conclusion}

Nonmatching quadratures may lead to nonconvergence of policy iteration.
In order to retain the fast superlinear convergence of policy iteration, matching quadrature should be used for all terms relating to the Hamiltonian.

\bibliographystyle{siamplain_NoURL}
\bibliography{refs}

\end{document}